\documentclass{article}

\usepackage{amsmath,amssymb}	

\usepackage{graphicx}		  
\usepackage[all]{xy}
\usepackage{tjm}

\newtheorem{theorem}{Theorem}[section]
\newtheorem{corollary}[theorem]{Corollary}
\newtheorem{proposition}[theorem]{Proposition}
\newtheorem{lemma}[theorem]{Lemma}

\newcommand{\bb}[1]{\mathbb{#1}}
\newcommand{\mcal}[1]{\mathcal{#1}}
\newcommand{\bbZ}{\bb{Z}}

\newcommand{\mrm}[1]{\mathrm{#1}}
\newcommand{\mbf}[1]{\mathbf{#1}}

\newcommand{\Spec}{\mrm{Spec}}

\newcommand{\mfr}[1]{\mathfrak{#1}}

\newcommand{\vare}{\varepsilon}

\title[Exponential map for Pic]{Notes on a $p$-adic exponential map for the Picard group}

\author{Wataru Kai}

\affils{Universit\"{a}t Essen}



\authorname{Wataru Kai}
\address{Fakult\"{a}t f\"{u}r Mathematik,\\
Universit\"{a}t Essen,\\
Thea-Leymann-Str.\ 9, 45127 Essen, Germany.}
\email{kaiw@ms.u-tokyo.ac.jp}

\begin{document}

\maketitle

\begin{abstract}
For proper flat schemes over a complete discrete valuation ring of mixed characteristic,
we construct an isomorphism of certain subgroups
of the Picard group and the first cohomology group of the structure sheaf.
When the Picard scheme is available and smooth, it recovers the isomorphism coming from its formal completion.
A reinterpretation of an old theorem of Mattuck is given.
\end{abstract}

\section*{Introduction}

For a compact complex manifold $X$, one has a short exact sequence of abelian sheaves:
\[  0\to 2\pi i \bbZ \to \mcal{O}_X \xrightarrow{\exp }\mcal{O}_X^* \to 1  \]
which gives an identification of a subgroup of $\mrm{Pic}(X)$ with a complex torus $H^1(X,\mcal{O}_X)/\mrm{Im}\ (2\pi i\ H^1(X, \bbZ ))$.
For a proper smooth geometrically connected variety $X$ over $\bb{Q}_p$ we have an analogue of this: We know that $\mrm{Pic}(X)$ is the group of $\bb{Q}_p$-points of a group scheme $\mbf{Pic}_X $ whose neutral component $\mbf{Pic}^0_X$ is an abelian variety, say of dimension $g$; by a theorem of Mattuck \cite{Mattuck}, the group of $\bb{Q}_p$-points of an abelian variety over $\bb{Q}_p$ contains an open subgroup (with respect to the valuation topology) isomorphic to $\bbZ _p^{\oplus g}$. 
In fact, by his proof we can take a sufficiently small open $\bbZ _p$-submodule $I $ of $H^1(X,\mcal{O}_X)$ and a canonical open immersion $I \hookrightarrow \mrm{Pic}_X^0(\bb{Q}_p)$.
These examples realize on a concrete level the fact that $H^1(X,\mcal{O}_X)$ is the tangent space of the Picard functor.
If we consider $\bbZ _p^{\oplus g}$ as a $p$-adic version of complex tori, the two stories are parallel also in that the Picard group of a proper smooth variety contains a torus as an open subgroup.

In these notes we describe how to construct the latter local isomorphism of $H^1(X,\mcal{O}_X)$ and $\mrm{Pic}(X)$ in the $p$-adic context in such a way that makes this parallel clearer:
we directly use the exponential power series to construct it.
It works for an arbitrary complete discrete valuation ring $\Lambda $ faithfully flat over $\bb{Z}_p$ and a proper $\Lambda [1/p]$-scheme, once one chooses a proper flat model for the scheme.
A notable feature of this construction is that we do not need the Picard variety at all, which may not exist for a general proper scheme.

\paragraph{Structure of these notes.} 

The set-up and the construction are explained in \S \ref{Sec.2}.
We shall consider a complete discrete valuation ring $(\Lambda ,\mfr{m})$ faithfully flat over $\bb{Z}_p$ and a proper flat $\Lambda $-scheme $X$. We will construct an isomorphism between $\ker [\mrm{Pic}(X)\to \mrm{Pic}(X\otimes _{\Lambda }\Lambda /\mfr{m}^n)]$ and $\ker [H^1(X,\mcal{O}_X)\to H^1(X\otimes _{\Lambda }\Lambda /\mfr{m}^n)]$ for sufficiently large $n$.
A bound for $n$ is given in terms of cohomological information.
In the same range, these groups are shown to inject into the corresponding groups for the generic fiber.

When the Picard functor for the scheme is formally smooth along the zero section, another exponential map is defined via formal group theory. In \S \ref{Sec.3} we verify that this map can be recovered from our construction.
For this we use the functoriality of our maps with respect to flat scalar extensions $\Lambda \to S$ to $p$-adically complete local rings $S$ which are not necessarily discrete valuation rings.

In \S \ref{Sec.4} we give a reinterpretation of Mattuck's proof of the above-mentioned theorem using our log/exp map.

\section{The exponential map for the Picard group}\label{Sec.2}

We begin with simple facts.
Let $R$ be a $p$-adically complete northerian ring flat over $\bb{Z}_{p}$, and $I$ be a defining ideal that is principal and $e\ge 1$ be an integer such that $p\in I^e$ (for example one can always take $I=(p)$ and $e=1$).

Consider the power series $\exp (T)=\sum _{i\ge 0} \frac{T^i}{i!}$ and
$\log (1-T)=\sum _{i\ge 1} \frac{T^i}{i}$.

\begin{lemma}\label{Lem.series}
Let $n$ be an integer satisfying $n > \frac{e}{p-1}$. Then for every element $x\in I^n$ and $i\ge 1$, the element $\frac{x^i}{i!}\in I^n$ makes sense and the infinite sum $\exp (x)= \sum _{i\ge 0} \frac{x^i}{i!}$ converges in $R$. Moreover we have $\exp (x)\equiv 1+x \mod (I^{n+1})$.

Similarly, for such an $x$ the sum $\log (1-x)=\sum _{i\ge 1}\frac{x^i}{i}$ converges in $R$ and we have $\log (1-x)\equiv x \mod (I^{n+1})$.
\end{lemma}
\begin{proof}
If $p^{n(i)}$ is the highest power of $p$ that divides $i!$, we have $n(i)= [\frac{i}{p}]+[\frac{i}{p^2}]+ \cdots < \frac{i}{p-1}$ (unless $i=0$). The assertion follows from this.
\end{proof}

\begin{corollary}
Let $n>\frac{e}{p-1}$.
We have isomorphisms of abelian groups (the latter of which is a subgroup of $R^*$)
\[ I^n \underset{\log }{\overset{\exp }{\rightleftarrows }} 1+I^n.  \]
\end{corollary}

Let $\mfr{X}$ be a formal scheme that is adic over $R$. This means that the formal affine open subsets of $\mfr{X}$ have coordinate rings which are $I$-adically topologized. The above lemmas hold for formal affine open sets of $\mfr{X}$ (because the hypotheses are again true for them), so we have isomorphisms of abelian groups on $\mfr{X}$ for $n>\frac{e}{p-1}$
\[ I^n\mcal{O}_{\mfr{X}} \underset{\log }{\overset{\exp }{\rightleftarrows }} 1+I^n\mcal{O}_{\mfr{X}}    \]
which gives an isomorphism $H^1(\mfr{X},I^n\mcal{O}_{\mfr{X}})\cong H^1(\mfr{X}, 1+I^n\mcal{O}_{\mfr{X}})$. We want to relate these groups with the honest cohomology groups $H^1(\mfr{X},\mcal{O}_{\mfr{X}})$ and $H^1(\mfr{X}, \mcal{O}^*_{\mfr{X}})$,
when the base ring is a complete discrete valuation ring.

\subsection{Set-up and the isomorphism.}\label{Sec.2.1}
Let $\Lambda $ be a complete discrete valuation ring flat over $\bb{Z}_{p}$ and let $\mfr{m}$ be its maximal ideal.
Set the integer $e$ by the relation $(p)=\mfr{m}^e$.
Let $X$ be a proper flat $\Lambda $-scheme. 

By \cite[III${}_2$ Th.7.7.6]{EGA}, there is a $\Lambda $-module $M_X$ of finite type such that for all $\Lambda $-algebras $S$, we have functorial isomorphisms
\[ H^0(X\otimes _\Lambda S,\mcal{O}_{X\otimes _\Lambda S}) \cong \mrm{Hom}_{\Lambda \text{-mod}}(M_X,S). \]
Denote its torsion part by $M_X^{\mrm{tors}}$ and let $t(X)\ge 0$ be the integer such that $\mfr{m}^{t(X)}=\mrm{Ann}_{\Lambda }(M_X^{\mrm{tors}})$. 

\begin{lemma}\label{KeepTheNotation}
Keep the notation above.
Let $S$ be a flat noetherian $\Lambda $-algebra.
If $n>t(X)$, then the group $H^0(X\otimes _{\Lambda }S/\mfr{m}^nS,\mcal{O}_{X\otimes _{\Lambda }S/\mfr{m}^nS})$ is generated by the images of $H^0(X\otimes _\Lambda S,\mcal{O}_{X\otimes _\Lambda S})$ and $H^0(X\otimes _\Lambda S/\mfr{m}^nS ,\mfr{m}^{n-t(X)}\mcal{O}_{X\otimes _\Lambda S/\mfr{m}^nS })$.
\end{lemma}
\begin{proof}
We know $H^0(X\otimes _\Lambda S/\mfr{m}^nS,\mcal{O}_{X\otimes _\Lambda S/\mfr{m}^nS})\cong \mrm{Hom}_{\Lambda \text{-mod}}(M_X,S/\mfr{m}^{m}S)$.
From the facts that the $p$-primary torsion part of $M_X$ is annihilated by $\mfr{m}^{t(X)}$, our $S$ is flat over $\Lambda $, and $n>t(X)$, it follows that any $\varphi \in \mrm{Hom}_{\Lambda \text{-mod}}(M_X,S/\mfr{m}^nS)$ must map the torsion part of $M_X$ into $\mfr{m}^{n-t(X)}S/\mfr{m}^nS$.
Let $\varphi ^{\mrm{tors}}\colon M_X^{\mrm{tors}}\to \mfr{m}^{n-t(X)}S/\mfr{m}^nS $ the induced map.

Put $M_X^{\mrm{free}}:= M_X/M_X^{\mrm{tors}}$ which is a free $\Lambda $-module. We know $M_X$ admits a direct sum decomposition $M_X\cong M_X^{\mrm{tors}}\oplus M_X^{\mrm{free}}$. Under this chosen isomorphism, put $\varphi ^{\mrm{free}}\colon M_X^{\mrm{free}}\to S/\mfr{m}^nS$ be the induced map by $\varphi $. 
Thus our $\varphi $ is written as $\varphi = \varphi ^{\mrm{tors}} \oplus \varphi ^{\mrm{free}}$. 
Since $M_X^{\mrm{free}}$ is free, the map $\varphi ^{\mrm{free}}$ lifts to a map $M_X^{\mrm{free}}\to S$. This competes the proof.
\end{proof}

\begin{corollary}\label{Corollary}
In addition to the hypotheses of Lemma \ref{KeepTheNotation}, assume that $p$ belongs to the Jacobson radical of $S$ (which is the case e.g. when $S$ is $p$-adically complete).
Then the group $H^0(X\otimes _{\Lambda }S/\mfr{m}^nS, \mcal{O}^*_{X\otimes _{\Lambda }S/\mfr{m}^nS})$ is generated by the images of 
$H^0(X\otimes _\Lambda S,\mcal{O}^*_{X\otimes _\Lambda S})$ and $H^0(X\otimes _\Lambda S/\mfr{m}^n, 1+\mfr{m}^{n-t(X)} \mcal{O}_{X\otimes _\Lambda S/\mfr{m}^n}S)$.
\end{corollary}
\begin{proof}
Take any $a\in H^0(X\otimes _{\Lambda }S/\mfr{m}^nS, \mcal{O}^*_{X\otimes _{\Lambda }S/\mfr{m}^nS})$. 
By Lemma \ref{KeepTheNotation}, it is written as
\[ a=a_0+a_1  \] 
where $a_0 $ comes from the global section on $X\otimes _\Lambda S$ and $a_1$ is a section of the subsheaf $\mfr{m}^{n-t(X)}\mcal{O}_{X\otimes _\Lambda S/\mfr{m}^nS}$ of $\mcal{O}_{X\otimes _\Lambda S/\mfr{m}^nS}$. 
Since $H^0(X\otimes _\Lambda S,\mcal{O}_{X\otimes _\Lambda S})$ is finite and flat over $S$, its Jacobson radical contains the image of $\mfr{m}$. It follows that $a_0$ is necessarily invertible already on $X\otimes _\Lambda S$.

Hence we can write the right hand side as $a_0(1+\frac{a_1}{a_0})$.
This completes the proof.
\end{proof}

\newcommand{\hattimes}{\hat{\otimes }}

Now, for an abelian group valued presheaf $F$ on the category of $\Lambda $-schemes, we write $F(X)^{(n)}:=\ker [F(X)\to F(X\otimes _{\Lambda }\Lambda /\mfr{m}^n)]$.

\begin{proposition}
Let $X$ be as at the beginning of this section \S \ref{Sec.2.1} and $\mfr{X}$ be its $p$-adic formal completion. Let $S$ be a flat noetherian $\Lambda $-algebra which is $p$-adically complete.
Then the isomorphisms
$H^1(\mfr{X}\hattimes _\Lambda S ,\mfr{m}^n\mcal{O}_{\mfr{X}\hattimes _\Lambda S})\cong H^1(\mfr{X}\hattimes _\Lambda S, 1+\mfr{m}^n\mcal{O}_{\mfr{X}\hattimes _\Lambda S})  $
gotten in Lemma \ref{Lem.series} induce isomorphisms of the groups
$ 
H^1(X\otimes _\Lambda S,\mcal{O}_{X\otimes _\Lambda S})^{(n)}$
and
$ 
\mrm{Pic}(X\otimes _\Lambda S)^{(n)}$
for $n>t(X)+\frac{e}{p-1}$.
\end{proposition}
\begin{proof}
By the short exact sequence $0\to \mfr{m}^n\mcal{O}_{\mfr{X}\hattimes _\Lambda S}\to \mcal{O}_{\mfr{X}\hattimes _\Lambda S}\to \mcal{O}_{\mfr{X}\hattimes _\Lambda S/\mfr{m}^nS}\to 0$ and its multiplicative counterpart, we have two parallel exact sequences (cohomology groups are taken for the formal scheme $\mfr{X}\hattimes _\Lambda S$):
\footnotesize
\[  \xymatrix{
H^0(\mcal{O}_{\mfr{X}\hattimes _\Lambda S})\xrightarrow{r_a} H^0(\mcal{O}_{\mfr{X}\hattimes _\Lambda S/\mfr{m}^nS})\ar[r]^(0.65){\partial _a}
&H^1(\mfr{m}^n\mcal{O}_{\mfr{X}\hattimes _\Lambda S})\ar[r] 
&H^1(\mcal{O}_{\mfr{X}\hattimes _\Lambda S})^{(n)} \to 0
\\%
H^0(\mcal{O}^*_{\mfr{X}\hattimes _\Lambda S})\xrightarrow{r_m} H^0(\mcal{O}^*_{\mfr{X}\hattimes _\Lambda S/\mfr{m}^nS})\ar[r]^(0.6){\partial _m}
&H^1(1+\mfr{m}^n\mcal{O}_{\mfr{X}\hattimes _\Lambda S})\ar[r] \ar@{-}[u]_{\cong }
&\mrm{Pic}({\mfr{X}\hattimes _\Lambda S})^{(n)} \to 0
}  \]
\normalsize
where the vertical isomorphism coming from $\exp$/$\log $ is available whenever $n> \frac{e}{p-1}$. 
We know all the terms but $H^1(1+\mfr{m}^n\mcal{O}_{\mfr{X}\hattimes _\Lambda S})$ have counterparts for the usual scheme $X\otimes _\Lambda S$ and are respectively isomorphic to them
(\cite[III${}_1$ Prop.5.1.2]{EGA} for coherent cohomology,
\cite[III${}_1$ Cor.5.1.6, I (10.11.7.1), $0_{\mbf{I}}$ \S (5.4.3)]{EGA} for $\mrm{Pic}$).

We claim that if $n>t(X)+\frac{e}{p-1}$, the images of $\partial _a$ and $\partial _m$ coincide under the indicated isomorphism.

By Lemma \ref{KeepTheNotation} and Corollary \ref{Corollary}, the maps $\iota _a\colon H^0(\mfr{m}^{n-t(X)}\mcal{O}_{X\otimes _\Lambda S/\mfr{m}^nS}) \to \mrm{cok}(r_a)$ and $\iota _m\colon H^0(1+\mfr{m}^{n-t(X)}\mcal{O}_{X\otimes _\Lambda S/\mfr{m}^nS})\to \mrm{cok}(r_m)$ are surjective.
Therefore in order to verify our claim, it suffices to note that we have a commutative diagram with $\exp $/$\log $ vertical isomorphisms (which one checks easily by \v{C}ech description)
\[  \xymatrix{
H^0(\mfr{m}^{n-t(X)}\mcal{O}_{X\otimes _\Lambda S/\mfr{m}^nS}) \ar@{-}[d]_{\cong } \ar[r]^(0.6){\partial _a\circ \iota _a}& H^1(\mfr{m}^n\mcal{O}_{X\otimes _\Lambda S})\ar@{-}[d]^{\cong }
\\%
H^0(1+\mfr{m}^{n-t(X)}\mcal{O}_{X\otimes _\Lambda S/\mfr{m}^nS}) \ar[r]_(0.6){\partial _m\circ \iota _m}& H^1(1+\mfr{m}^n\mcal{O}_{\mfr{X}\hattimes _\Lambda S}).
}   \]
Here, the left hand vertical isomorphism exists because $n-t(X)>\frac{e}{p-1}$.
\end{proof}

We denote by $\exp _X(S)$ the isomorphism thus gotten:
\[ H^1(X\otimes _\Lambda S,\mcal{O}_{X\otimes _\Lambda S})^{(n)} \xrightarrow{\exp _X(S)}\mrm{Pic }(X\times _\Lambda S)^{(n)} .  \]
From the construction, it is functorial for flat noetherian $\Lambda $-algebras $S$ which are $p$-adically complete.

Next we show both hand sides injects into the corresponding groups for the generic fiber in pleasant circumstances.

\subsection{Injectivity.}

\begin{proposition}\label{PropKeep}
Keep the notation in \S \ref{Sec.2.1}.
For any flat noetherian $\Lambda $-algebra $S$, the restriction map
\[  H^1(X\otimes _\Lambda S,\mcal{O}_{X\otimes _\Lambda S})^{(n)} \to H^1(X\otimes _\Lambda S[ {1}/{p}] ,\mcal{O}_{X\otimes _\Lambda S[{1}/{p}]})   \]
is injective whenever $n> t(X)$.
\end{proposition}
\begin{proof}
Let us recall \cite[III${}_2$ Rem.6.10.6]{EGA} that there is a bounded complex
\[  0\to L^0\xrightarrow{d^0} L^1\to \cdots    \]
of finitely generated free $\Lambda $-modules such that there are functorial quasi-isomorphisms for all $\Lambda $-algebras $S$:
\[ L^\bullet \otimes _\Lambda S\xrightarrow{\sim } R\Gamma (X\otimes _\Lambda S, \mcal{O}_{X\otimes _\Lambda S}).   \]
(In fact, our $M_X$ in \S \ref{Sec.2.1} is obtained as the cokernel of the $\Lambda $-linear dual of the map $L^0\to L^1$.)
Let $\pi \in \mfr{m}$ be a uniformizer of $\Lambda $. By Elementary Divisors Theorem in linear algebra, we may assume, up to (zig-zag of) quasi-isomorphism and coordinate change of $L^0$ and $L^1$, that $d^0$ is defined by a diagonal matrix with entries $(\pi ^{e_1},\dots ,\pi ^{e_r}, 0, \dots )$ where $1\le e_1\le \dots \le e_r=t(X)$.

Therefore $H^1(X\otimes _\Lambda S,\mcal{O}_{X\otimes _\Lambda S})$ is a submodule of (with some $a\ge 0$)
\[ ( \bigoplus _{k=1}^r S/\pi ^{e_k}S ) \oplus S^a \] 
and $H^1(X\otimes _\Lambda S/\pi ^nS,\mcal{O}_{X\otimes _\Lambda S/\pi ^nS})$ 
is a submodule of $( \bigoplus _{k=1}^r S/\pi ^{e_k}S )\oplus (S/\pi ^nS)^a $ because $n> t(X)$. It follows that $H^1(X\otimes _\Lambda S,\mcal{O}_{X\otimes _\Lambda S})^{(n)}$
is a submodule of $(\pi ^n S)^a$ which is $\pi $-torsion free. Hence the restriction map $H^1(X\otimes _\Lambda S,\mcal{O}_{X\otimes _\Lambda S})^{(n)}$ $\to $ $H^1(X\otimes _\Lambda S[{1}/{p}],\mcal{O}_{X\otimes _\Lambda S[{1}/{p}]})$ is injective.
\end{proof}

\begin{proposition}
In the notation of Proposition \ref{PropKeep}, assume moreover that $S$ is a discrete valuation ring. Then the restriction map
\[ \mrm{Pic}(X\otimes _\Lambda S)^{(n)}\to \mrm{Pic}(X\otimes _\Lambda S[1/p]) \] 
is injective for all $n> t(X)$.
\end{proposition}
\begin{proof}
This is \cite[Lem.6.4.4]{Raynaud}.
\end{proof}

To summarize, for $p$-adically complete flat noetherian $\Lambda $-algebras $S$, we have established diagrams 
\begin{equation}\label{established}
\xymatrix{
H^1(X\otimes _\Lambda S[1/p],\mcal{O}_{X\otimes _\Lambda S[1/p]}) 
& \mrm{Pic}(X\otimes _\Lambda S[1/p])
\\%
H^1(X\otimes _\Lambda S,\mcal{O}_{X\otimes _\Lambda S})^{(n)} \ar@{-}[r]_(0.55){\cong } \ar@{^{(}->}[u]
& \mrm{Pic}(X\otimes _\Lambda S)^{(n)}
\ar[u]
} \end{equation}
for $n>t(X)+\frac{e}{p-1}$, functorial in $S$.
Moreover, if $S$ is a discrete valuation ring, the right hand vertical map is injective.

\subsection{Independence of the integral model.}
Suppose a proper scheme $X^o$ over $\Lambda [1/p]$ is given. Choose a proper flat $\Lambda $-scheme $X$ which contains $X^o$ as the generic fiber. One way to see the diagram \eqref{established} is that the subgroups $(-)^{(n)}$ topologize the groups $H^1(X^o,\mcal{O}_{X^o})$ and $\mrm{Pic}(X^o)$, and both hand sides are locally isomorphic for these topologies. 
One can ask if these topologies and the isomorphism are independent of the model $X$ we choose. 
(Of course, the subgroups $H^1(X,\mcal{O}_X)^{(n)}$ {\em are} open subgroups defining the usual topology on $H^1(X^o,\mcal{O}_{X^o})$ as a finite dimensional $\Lambda [1/p]$-vector space.)

The answer is that they are independent up to passage to smaller domain of definition of the local isomorphism. 
This follows from the next observation whose verification is immediate.

\begin{proposition}
Suppose two schemes $X\to \mrm{Spec}(\Lambda )$ and $X'\to \mrm{Spec}(\Lambda ')$ as in \S \ref{Sec.2.1} and a commutative diagram 
\[ \xymatrix{
X' \ar[d]\ar[r]& X \ar[d] \\%
\mrm{Spec}(\Lambda ')\ar[r]& \mrm{Spec}(\Lambda )
} \]
are given. Then the resulting diagram of abelian groups
\[ \xymatrix{
H^1(X',\mcal{O}_{X'})^{(n')}\ar[rr]_{\cong }^{\exp _{X'}(\Lambda ')}&& \mrm{Pic}(X')^{(n')} \\%
H^1(X,\mcal{O}_{X})^{(n)} \ar[u]
\ar[rr]^{\cong }_{\exp _{X}(\Lambda )}&& \mrm{Pic}(X)^{(n)} \ar[u]
}  \]
is commutative whenever the maps involved are defined.
\end{proposition}

Now, suppose we have two proper flat models $X_1,X_2$ of $X^o$. Let $X_3$ be the scheme theoretic closure of the diagonal $X^o\hookrightarrow X_1\times _{\Lambda }X_2$. It is proper and flat over $\Lambda $. We apply the previous proposition to the $\Lambda $-maps $X_3\to X_1$ and $X_3\to X_2$ to get a commutative diagram
\[ \xymatrix{
{}&H^1(X_1,\mcal{O}_{X_1})^{(n_1)} 
\ar[d]
\ar@{^{(}->}[dl] \ar@{-}[r]_{\cong }
&\mrm{Pic}(X_1)^{(n_1)}\ar@{^{(}->}[rd] \ar[d]
\\%
H^1(X^o,\mcal{O}_{X^o})
&H^1(X_3,\mcal{O}_{X_3})^{(n_3)}\ar@{_{(}->}[l]\ar@{-}[r]_{\cong }
&\mrm{Pic}(X_3)^{(i_3)}\ar@{^{(}->}[r]
&\mrm{Pic}(X^o)
\\%
{}&H^1(X_2,\mcal{O}_{X_2})^{(n_2)}
\ar[u]
\ar@{_{(}->}[ul]\ar@{-}[r]_{\cong }
&\mrm{Pic}(X_2)^{(n_2)}\ar@{^{(}->}[ru] \ar[u]
} \]
with $n_1,n_2,n_3$ large enough.
Recall that the indicated subgroups of $H^1(X^o,\mcal{O}_{X^o})$ are indeed open subgroups in the usual topology as a $\Lambda [1/p]$-vector space.
So we get the above-mentioned independence.

\section{Comparison with the formal exponential map}\label{Sec.3}

Keep the notation at the beginning of \S \ref{Sec.2.1}.
In this section we {\em assume} that the Picard functor for $X/\Lambda $ is formally smooth along the zero section.
This means the following: Let $(\mrm{Art}/\Lambda )$ be the category of local Artinian $\Lambda $-algebras $(S,\mfr{m}_S)$ equipped with an isomorphism $S/\mfr{m}_S \xrightarrow{\cong } \Lambda /\mfr{m}$. Then the covariant functor
\[  \mrm{Pic}^{\wedge }_{X/\Lambda }\colon
(\mrm{Art}/\Lambda )\ni S \mapsto
\ker [\mrm{Pic}(X\otimes _\Lambda S) \to \mrm{Pic}(X\otimes _\Lambda \Lambda /\mfr{m}) ]  \]
is co-represented by a formal power series ring $\Lambda [[T_1,\dots ,T_g]]$
(a pro-object in $(\mrm{Art}/\Lambda )$) 
as a set-valued functor.
As usual, such a functor is called a formal group over $\Lambda $.

Every finitely generated free $\Lambda $-module $L$ determines a formal group $L^+\colon $ $(S,\mfr{m}_S)\mapsto L\otimes _\Lambda \mfr{m}_S$.
For a formal group $F$ and an ideal $I\subset \mfr{m}_S$, the symbol $F(I)$ denotes $\ker [F(S)\to F(S/I)]$.

By the identification of the tangent space of $\mrm{Pic}^{\wedge }_{X/\Lambda }$ with $H^1(X,\mcal{O}_{X})$, we find that this cohomology group is torsion-free under our assumption. Since $\Lambda [1/p]$ is a characteristic $0$ field, the theory of formal group laws gives an exponential isomorphism
\[ \mrm{Exp}_X\colon H^1(X,\mcal{O}_{X})^+ \xrightarrow{\cong } \mrm{Pic}^{\wedge }_{X/\Lambda }   \]
{\em after} scalar extension to $\Lambda [1/p]$. It is the {\em unique} homomorphism of formal group laws inducing the given identification of the tangent spaces.
From this, for a large enough power $\mfr{m}^n$ we get an isomorphism of abelian groups
\[  \mrm{Exp}_X^{(n)}\colon
H^1(X,\mcal{O}_{X})\otimes _\Lambda \mfr{m}^n \xrightarrow{\cong } \mrm{Pic}^{\wedge }_{X/\Lambda }(\mfr{m}^n) \overset{\text{by def}}{=} \mrm{Pic}(X)^{(n)} . \]
Note that the general theory does not tell us how large $n$ has to be.

We are going to show that the map $\mrm{Exp}_X$ can be constructed from our $\exp _X$, see Theorem \ref{ExpRecover}. 
Then we find {\it a posteriori} 
a bound of $n$ for the existence of $\mrm{Exp}_X^{(n)}$: it exists whenever $n>t(X)+\frac{e}{p-1}$.
Notice that from the proof of Proposition \ref{PropKeep} we have $H^1(X,\mcal{O}_X)\otimes _\Lambda \mfr{m}^n=\mfr{m}^nH^1(X,\mcal{O}_X)=H^1(X,\mcal{O}_X)^{(n)}$ for $n\ge t(X)$.

\subsection{Extracting power series.}\label{Extracting}

Choose an isomorphism 
\begin{equation}\label{ChosenIso} 
\mrm{Pic}^{\wedge }_{X/\Lambda }\cong \mrm{Hom}_{\mrm{Art}/\Lambda }(\Lambda [[T_1,\dots ,T_g]],-) \end{equation} 
of set-valued functors.
The abelian group law on $\mrm{Pic}^{\wedge }_{X/\Lambda }$ is translated into a map $\Lambda [[T_1,\dots, T_g]]\to \Lambda [[T_1,\dots ,T_g ]]\hattimes \Lambda [[T_1,\dots ,T_g ]] $.
Denote by $\tau _j(T_1\otimes 1, \dots ,1\otimes T_g)$ the image of $T_j$.
The isomorphism gives an isomorphism $H^1(X,\mcal{O}_X)\cong \Lambda ^{\oplus g}$ through the identification of tangent spaces.
Denote by $e_1,\dots ,e_g\in H^1(X,\mcal{O}_X)$ the basis thus given.
Choose a uniformizer $\pi $ of $\Lambda $. Let $n> t(X)+\frac{e}{p-1}$.

Let $t_1,\dots ,t_g$ be another set of indeterminates. The ring $\Lambda [[t_1,\dots ,t_g]]$ is $p$-adically complete, so we can apply the results of \S \ref{Sec.2} to get a map
\begin{multline*}  
\exp _{X}(\Lambda [[\underline{t}]]) \colon H^1(X,\mcal{O}_X)^{(n)}\otimes _\Lambda \Lambda [[t_1,\dots ,t_g]] \to
\mrm{Pic}(X\otimes _\Lambda \Lambda [[t_1,\dots ,t_g]])^{(n)}
\\%
\cong \ker \biggl(  \begin{array}{l}
\mrm{Hom}_{\mrm{(Art)}/\Lambda }(\Lambda [[T_1,\dots ,T_g]], \Lambda [[t_1,\dots ,t_g]]) \\
\to \mrm{Hom}_{\mrm{(Art)}/\Lambda }(\Lambda [[T_1,\dots ,T_g]],\Lambda /\mfr{m}^n [[t_1,\dots ,t_g]]) \end{array}
\biggr) . \end{multline*}
The map is identified with a map, via our bases:
\[ \exp _{X}(\Lambda [[\underline{t}]]) 
\colon
\bigoplus _{j=1}^g e_j \mfr{m}^n[[t_1,\dots ,t_g]] \to \prod _{j=1}^g \mfr{m}^n[[t_1,\dots ,t_g]].   \]

For each $j$, we consider the $j$-th entry of the image of the element $\sum _{k=1}^ge_k \pi ^nt_k $ which we write as $\pi ^n f_j(t_1,\dots ,t_g)$ with $f_j\in \Lambda [[t_1,\dots ,t_g]]$.

From the functoriality it follows that given any local complete noetherian faithfully flat $\Lambda $-algebra $(S,\mfr{m}_S)$ and elements $x_1,\dots ,x_g\in \mfr{m}_S$ (so $f_j(x_1,\dots ,x_g)\in S$ make sense), our map
\[ \exp _X(S)\colon \mfr{m}^nS\ H^1(X,\mcal{O}_X) \to \mrm{Pic}(X\otimes _\Lambda S)^{(n)}\cong \prod _{j=1}^g \mfr{m}^nS   \]
sends the element $\sum _{j=1}^g \pi ^nx_j e_j $ to $(\pi ^n f_j(x_1,\dots ,x_g))_j$.
The fact that $\exp _X(-)$ is a functorial homomorphism implies the identity
\[ \pi ^n f_j(x_1+y_1, \dots ,x_g+y_g)
= \tau _j (\pi ^n f_1(x_1,\dots ,x_g), \dots ,
\pi ^nf_g(y_1,\dots ,y_g)). \]
This means that the power series $\pi ^nf_j $ $(j=1 ,\dots ,g)$ constitute a homomorphism of formal group laws 
$H^1(X,\mcal{O}_X)^+\to \mrm{Pic}^{\wedge }_{X/\Lambda }$
over $\Lambda $.

\subsection{The comparison.}

In the rest of the section we prove:

\begin{theorem}\label{ExpRecover}
Let $(E_j(T_1,\dots ,T_g))_{j=1}^g$ be the power series with $\Lambda [1/p]$-coefficients defining $\mrm{Exp}_X$ at the beginning of this \S \ref{Sec.3}.
Let $(\pi ^n f_j )_j$ be the power series gotten in \S \ref{Extracting}.
Then we have the identity 
\[ E_j(\pi ^nT_1,\dots ,\pi ^nT_g)=\pi ^n f_j(T_1,\dots ,T_g). \]
\end{theorem}

By the uniqueness of $E_j$, it suffices to show that the degree $1$ part of $\pi ^nf_j$ is $\pi ^nT_j$.
We are going to exploit the following framework to show this.
An element in $\Lambda $ is determined by its mod $\pi ^\nu $ classes for all $\nu \ge 1$. Fix a $\nu \ge 1$. Choose $\mu $ very large.
The composite map
\[  \pi ^\mu H^1(X,\mcal{O}_X) \subset
H^1(X,\mcal{O}_X)\xrightarrow{(\pi ^n f_j)_j } \mrm{Pic}(X)^{(n)}  \]
which lands in $\mrm{Pic}(X)^{(n+\mu )}$,
followed by
\[ \mrm{Pic}(X)^{(n+\mu )}\twoheadrightarrow  
\frac{\mrm{Pic}(X)^{(n+\mu )}}{
\mrm{Pic}(X)^{(n+\mu +\nu )}  } \]
will depend only on the coefficients of degree $1$ of $\pi ^nf_j$ if $\mu $ is sufficiently large (depending on $\nu $ and coefficients of $\pi ^nf_j$), because then the higher terms of $\pi ^nf_j$ will have very high divisibility by $\pi $ compared with the degree $1$ terms.
Let us call this map $\Theta _{\mu }$: 
\[ \Theta _\mu \colon \pi ^\mu H^1(X,\mcal{O}_X)\to \frac{\mrm{Pic}(X)^{(n+\mu )}}{
\mrm{Pic}(X)^{(n+\mu +\nu )}  }.\]

We will show that this map is equal to the map
\[ (\mfr{m}^\mu )^g \xrightarrow{\times \pi ^n} (\mfr{m}^{n+ \mu })^g \overset{\mrm{can.}}{\twoheadrightarrow }(\mfr{m}^{n+ \mu }/\mfr{m}^{n+ \mu +\nu })^g \]
under identifications $\pi ^\mu H^1(X,\mcal{O}_X)=\bigoplus _{j=1}^g \mfr{m}^\mu e_j$ and 
$ \frac{\mrm{Pic}(X)^{(n+\mu )}}{
\mrm{Pic}(X)^{(n+\mu +\nu )}  } 
\cong 
\biggl( \frac{\mfr{m}^{n+\mu }}{\mfr{m}^{n+\mu +\nu } } \biggr) ^g.$
Then we will be able to conclude that the degree $1$ part of $\pi ^nf_j $ is $\pi ^nT_j$ mod $\mfr{m}^\nu $.

In order to carry it out, recall that the canonical identification of the tangent space of $\mrm{Pic}^{\wedge }_{X/\Lambda }$ with $H^1(X,\mcal{O}_X)$
(which is behind our choice of the basis $(e_j)_j$)
is obtained from the split short exact sequence
\[ 0\to \varepsilon \mcal{O}_{X[\vare ]/(\vare ^2)} \xrightarrow{+1} 
\mcal{O}^*_{X[\vare ]/(\vare ^2)}\to 
\mcal{O}^*_{X}
\to 1 \]
by taking $H^1(X[\vare ]/(\vare ^2),-)$.
Consider the next commutative diagram for $\mu \gg 0$:
\footnotesize
\[  \xymatrix{
0\to \varepsilon \mcal{O}_{X[\vare ]/(\vare ^2)} \ar[r]^{+1}\ar[d] &
\mcal{O}^*_{X[\vare ]/(\vare ^2)} \ar[r]
\ar[d] &
\ar[d] \ \mcal{O}^*_{X} 
\to 1
\\%
0\to \pi ^{n+\mu }\mcal{O}_{X\otimes _\Lambda \Lambda /\mfr{m}^{n+\mu +\nu } }  \ar[r]^{+1 }\ar@{-->}[dr]_{\exp }^{\cong } &
\mcal{O}^*_{X\otimes _\Lambda \Lambda /\mfr{m}^{n+\mu +\nu }} \ar[r]&
\mcal{O}^*_{X\otimes _\Lambda \Lambda /\mfr{m}^{n+ \mu } } \to 1
\\%
&1+\mfr{m}^{n+\mu }\mcal{O}_{X\otimes _\Lambda \Lambda /\mfr{m}^{n+\mu +\nu }} \ar@{}|{\cup }[u]
} \]
\normalsize
where the left and middle downward maps are given by $\vare \mapsto \pi ^{n+\mu } $.
It gives a commutative diagram
\small
\[  \xymatrix{
H^1(X,\mcal{O}_X) \ar@{-}[r]^{\cong }_{\mrm{can.}}
\ar[d]_{\times \pi ^{n+\mu }}
&T_0\mrm{Pic}^{\wedge }_{X/\Lambda }\ar[d]
\\%
H^1(X,\mfr{m}^{n+\mu }\mcal{O}_{X\otimes _\Lambda \Lambda /\mfr{m}^{n+\mu +\nu }})
\ar@{->>}[r] 
&\ker \biggl[ {\begin{array}{c} 
\mrm{Pic}(X\otimes _\Lambda \Lambda /\mfr{m}^{n+\mu +\nu }) \\ \downarrow \\ \mrm{Pic}(X\otimes _\Lambda \Lambda /\mfr{m}^{n+\mu })  \end{array}} \biggr]
_{\cong \biggl( \frac{\pi ^{n+\mu } \Lambda }{\pi ^{n+\mu +\nu }\Lambda } \biggr) ^g }
}  \]
\normalsize
where bijection at the lower right corner is the one that comes from the isomorphism \eqref{ChosenIso} we chose.
Thus tracking the upper path, we find that the basis $(e_j)_j$ of the upper left corner maps to the standard basis of the lower right corner.
On the other hand, by the construction of our $\exp _X$, we may add the next flagment to the bottom of this diagram:
\[  \xymatrix{ 
{}&{}&& {}
\\%
\pi ^\mu H^1(X,\mcal{O}_X)\ar[r]^{\times \pi ^n}
&\frac{H^1(X,\mfr{m}^{n+\mu } \mcal{O}_X)}{H^1(X,\mfr{m}^{n+\mu +\nu } \mcal{O}_X)}
\ar[rr]_{\exp _X}^{\cong } \ar[u]
&& \frac{\mrm{Pic}(X)^{(n+ \mu )}}{\mrm{Pic}(X)^{(n+\mu +\nu )}}
\ar[u]_{\cong }
\\%
&\bigoplus _{j=1}^g  \frac{\pi ^{n+\mu \Lambda }}{\pi ^{n+\mu +\nu }\Lambda }e_j \ar@{=}[u] &&
}  \]
The last horizontal composite computes
$\Theta _\mu $.
Tracking the bases, we find that the map $\Theta _\mu $ is described by the multiplication map by $\pi ^n$.

This completes the proof of the identity $E_j(\pi ^nT_1,\dots ,\pi ^nT_g)=\pi ^nf_j(T_1,\dots ,T_g)$.

\section{Recovering Mattuck's theorem}\label{Sec.4}

As before, let $(\Lambda ,\mfr{m})$ be a complete discrete valuation ring flat over $\bb{Z}_p$. Denote its fraction field by $k$.
Let $A$ be an abelian variety of dimension $g$ over $k$. Results applied to the dual abelian variety of $A$ implies that $A(k)$ contains a subgroup isomorphic to $\Lambda ^g$. 

This almost recovers Mattuck's theorem mentioned in Introduction (at least when the value group is {\em discrete}),
except
that it remains to be shown that our subgroup is open in the topology of $A(k)$ coming from that of $k$.
Here, for a locally finite-type $k$-scheme $X$, the set of its $k$-rational points $X(k)$ is topologized by the following procedure: we cover $X$ by affine open subsets $X'$, take a closed embedding $X'\hookrightarrow \bb{A}^N$ and give $X'(k)$ the subspace topology of $\bb{A}^N(k)=k^N$; the topology does not depend on choices.
It is often called the {\em valuation topology}.

One way to treat topology is to repeat his arguments partially.
We do this below.
As Mattuck (or anyone who tries to prove something equivalent e.g.\ through formal group theory) pays considerable efforts to handle the group structure, this process is much easier in our case.

\subsection{Reduction to Jacobian.}\label{ReductionToCurve}

Let $A^\vee $ be the dual abelian variety. By repeated hyperplane section we choose a closed subcurve $C$ of $A^\vee $, smooth and passing through a $k$-rational point.

By weak Lefschetz theorems, the map induced on the Picard varieties
$ A=\mbf{Pic}^0_{A^\vee }\hookrightarrow \mbf{Pic}^0_C $
is a closed immersion
and
$ H^1(A^\vee ,\mcal{O}_{A^\vee }) \hookrightarrow H^1(C,\mcal{O}_{C}) $
is injective.
In particular the objects attached to $A^\vee $ have the induced topology from those attached to $C$.
Therefore the problem is reduced to Jacobian varieties.

\subsection{Symmetric product and the Jacobian.}\label{NowLetC}

Now, let $C$ be a proper smooth curve over $k$ with genus $g$, having a $k$-rational point $x$. 
If $t$ is a local parameter of $C$ at $x$, an inverse map theorem tells us that the map $t\colon C(k)\to \bb{A}^1(k)$ defined around $x$ is an analytic isomorphism (i.e. the inverse map is given by a power series) of neighborhoods of $x$ and $0$.
Let $x_1,\dots ,x_g \in C(k)$ be distinct points, sufficiently general so that $\Gamma (C,\mcal{O}_C(x_1+\dots +x_g))\cong k$.

Let $C^{(g)}$ be the $g$-fold symmetric product of $C$. 
Choose any reference point $x_0\in C(k)$.
We know that the canonical map (depending on $x_0$ up to translation) $C^{(g)}\to \mbf{Pic}^0_{C}$ is birational. We further know that our $(x_1,\dots ,x_g)$ is inside the isomorphism locus under the current assumption, e.g.\ \cite[Th.5.1 or its proof]{Milne}.

\subsection{Choosing a model.}\label{Model}

By the desingularization of two-dimensional schemes, we can take a proper flat $\Lambda $-scheme $\mfr{C}$ containing $C$ as the generic fiber. Let $\mfr{x}_1,\dots ,\mfr{x}_g$ be the closure of the points $x_1,\dots ,x_g$ in $\mfr{C}$; they are isomorphic to $\Spec (\Lambda )$. Let $y_i$ be their respective closed point.
By blowing-up $\mfr{C}$ at some of these points, we may assume $y_1,\dots ,y_g$ are distinct.
Choose a defining equation $t_i$ of $\mfr{x}_i$ in $\Spec (\mcal{O}_{\mfr{C},y_i})$.
It gives a local parameter of $C$ at $x_i$ too.
For $1\le i\le g$ and $\vare \in \mfr{m}$, 
consider the subscheme 
\[ \mfr{x}_i(\vare ):=\Spec (\mcal{O}_{\mfr{C},y_i}/(t_i-\vare )) \]
which is isomorphic to $\Spec (\Lambda )$.
Denote its generic point by 
$ x_i(\vare ) \in C(k) $.

\subsection{First comparison.}

The Cartier divisor $\mfr{x}_i(\vare )-\mfr{x}_i$ on $\mfr{C}$ is defined by the local equation $(t_i-\vare )/t_i $ on $\Spec (\mcal{O}_{\mfr{C},y_i})$.
In particular if $\vare \in \mfr{m}^n$, its restriction to $\mfr{C}\otimes _{\Lambda }\Lambda /\mfr{m}^n$ is trivial.

We consider the next map of sets.
It maps $\prod _i \mfr{m}^n $ into $\mrm{Pic}(\mfr{C})^{(n)}$.
\[ \Phi \colon \prod _{i=1}^g \mfr{m} \to \mrm{Pic}(\mfr{C})\ ; \quad  
(\vare _i)_i 
\mapsto \mcal{O}_{\mfr{C}}( \sum _i (\mfr{x}_i(\vare )-\mfr{x}_i)) \]

The map $\Phi $ followed by the restriction map to $\mrm{Pic}(C)$ can be decomposed as 
\[ \prod _{i=1}^g\mfr{m}\to C^{(g)}\ ; \ 
(\vare _i)_{i=1}^g\mapsto (x(\vare _i))_{i} \]
followed by
\[ C^{(g)} \to \mbf{Pic}^0_C (k)\ ; \
(x'_i)_i \mapsto \mcal{O}_{C}(\sum _i (x_i' - x_i)). \]
Both of these maps are local homeomorphisms around the points we are looking at, by the facts in the previous \S \ref{NowLetC}.
So in order to achieve the purpose of this section, it suffices to show that the map $\Phi \colon \prod _{i=1}^g \mfr{m}\to \mrm{Pic}(\mfr{C})$ is a local homeomorphism around the origin.

Note that if $n$ is large enough, $\mrm{Pic}(\mfr{C})^{(n)}$ is a subset of $\mbf{Pic}^0_C(k)$ containing a subset
(the image of $\prod _i \mfr{m}^n$)
which is open in the valuation topology.
This already shows that the valuation topology is {\em finer} than the topology defined by $\mrm{Pic}(\mfr{C})^{(n)}$'s.
In order to prove that they are the same topology, we are going to look at the $\log $ isomorphism $\mrm{Pic}(\mfr{C})^{(n)}\xrightarrow{\cong } H^1(\mfr{C},\mcal{O}_{\mfr{C}})^{(n)}$.

\subsection{Choosing a coordinate system.}\label{Coordinate}

Put $D=\sum _{i=1}^g \mfr{x}_i$. It is an effective Cartier divisor on $\mfr{C}$ flat over $\Spec (\Lambda )$.
Twist by $D$ is denoted by $(-)(D)$.
As a $\Lambda $-scheme, $D$ is isomorphic to the disjoint union of $g$ copies of $\Spec (\Lambda )$, by the assumption made in \S \ref{Model} that $y_i$ are distinct.
On the scheme $\mfr{C}$ we have an exact sequence:
$ 0\to \mcal{O}_{\mfr{C}}\to 
\mcal{O}_{\mfr{C}}(D)\to \mcal{O}_{D}(D)
\to 0.  $
That gives an injection 
\[  \partial _D\colon \ (\Lambda ^g\cong \ )\  \bigoplus _i 
\biggl[ (\frac{1}{t_i}\mcal{O}_{\mfr{x}_i})/\mcal{O}_{\mfr{x}_i} \biggr] { \hookrightarrow } H^1(\mfr{C},\mcal{O}_{\mfr{C}}) \] 
with $p$-power torsion cokernel.

For an element $f$ in the total fraction ring of $\mcal{O}_{\mfr{C},y_i}$,
let us write $h^1(\mfr{x}_i,f ,\mcal{O}_{\mfr{C}} )\in H^1(\mfr{C},\mcal{O}_{\mfr{C}})$
for 
the element represented by the following \v{C}ech datum: the covering $\mfr{C}=\Spec (\mcal{O}_{\mfr{C},y_i})\cup (\mfr{C}\setminus |f|)$, where $|f|$ temporarily denotes the non-defined locus of $f$ in $\Spec (\mcal{O}_{\mfr{C},y_i})$, and the section 
$f$ on $\Spec (\mcal{O}_{\mfr{C},y_i})\setminus |f|$ of the structure sheaf.
Then for $\vare \in \Lambda $, 
one finds that
$h^1(\mfr{x}_i,\vare /t_i ,\mcal{O}_{\mfr{C}})$ is
the image of $\vare /t_i \in (\frac{1}{t_i}\mcal{O}_{\mfr{x}_i})/\mcal{O}_{\mfr{x}_i}$ under $\partial _D$.

Note that if $\vare $ belongs to a high enough power $\mfr{m}^n$ (suppose we fixed $n$) that annihilates $\mrm{coker}(\partial _D)$, then for any $m\ge 1$, the element $h^1(\mfr{x}_i,\vare /t_i^m ,\mcal{O}_{\mfr{C}})=\vare h^1(\mfr{x}_i,1 /t_i^m ,\mcal{O}_{\mfr{C}})\in H^1(\mfr{C},\mcal{O}_{\mfr{C}})$ belongs to the image of $\partial _D$, hence can be uniquely written as
\begin{equation}\label{uniquely} 
h^1(\mfr{x_i},\vare /t_i^m,\mcal{O}_{\mfr{C}})= \sum _{j=1}^g h^1(\mfr{x}_j, c(i,j, m,\vare )/t_j,\mcal{O}_{\mfr{C}})  \end{equation} 
with $c(i,j,m,\vare )\in \Lambda $. 
The element $c(i,j, m,\vare )$ depends linearly on $\vare $. If we write $c(i,j,m,\vare )= c(i,j,m)\vare $, then $c(i,j,m)\in k$ and has normalized valuation at least $-n$.

Note the tautological identity $c(i,j,1)=\delta _{ij}$ (Kronecker's delta).

\subsection{Rewriting $\Phi $.}\label{Rewriting}

Similarly let us denote by 
$h^1(\mfr{x}_i,f ,\mcal{O}_{\mfr{C}}^*)\in \mrm{Pic}(\mfr{C})$
the element 
determined by the \v{C}ech datum:
the covering $\mfr{C}=\Spec (\mcal{O}_{\mfr{C},y_i})\cup (\mfr{C}\setminus |f|)$, where $|f|$ now denotes the non-invertible locus of $f$ on $\Spec (\mcal{O}_{\mfr{C},y_i})$, and the section 
$f$ on $\Spec (\mcal{O}_{\mfr{C},y_i})\setminus |f|$ of the multiplicative group of invertible functions.
Then 
$h^1(\mfr{x}_i,1-\frac{\vare }{t_i},\mcal{O}_{\mfr{C}}^*)$ is equal to the class in $\mrm{Pic}(\mfr{C})$ represented by the Cartier divisor $\mfr{x}_i(\vare )-\mfr{x}_i$;
it follows 
\[ \Phi ((\vare _i)_i)
=\sum _{i=1}^g h^1(\mfr{x}_i,1-\frac{\vare _i}{t_i},\mcal{O}_{\mfr{C}}^*). \]

\subsection{End of comparison.}

Now let $n\ge 1$ be large enough so that the results in \S \ref{Sec.2} are available. Consider the composite
\[ \prod _{i=1}^g \mfr{m}^n \xrightarrow{\Phi } \mrm{Pic}(\mfr{C}  )^{(n)} \xrightarrow[\cong ]{\log }H^1(\mfr{C},\mcal{O}_{\mfr{C}})^{(n)}.  \]
Suppose our $n$ is large enough so the last group is contained in $\mrm{Im}(\partial _D)$.
By the previous \S \ref{Rewriting} and
the construction of $\log $,
we have 
\[ (\log \circ \Phi )\biggl( (\vare _i)_i \biggr) = \sum _{i=1}^g h^1(\mfr{x}_i,\log (1-\frac{\vare _i}{t_i}),\mcal{O}_{\mfr{C}})). \]
Here since $\log (1-\frac{\vare _i}{t_i})=\frac{\vare _i}{t_i}+\frac{1}{2}(\frac{\vare _i}{t_i})^2+\dotsb $ is a formal power series, this expression {\em a priori} makes sense only in the cohomology group of the $p$-adic formal completion of $\mfr{C}$, but by \cite[III${}_2$ Prop.5.1.2]{EGA} mentioned earlier it makes sense in the indicated cohomology group.

By the relation \eqref{uniquely} in \S \ref{Coordinate}, the right hand side can be written as:
\[ \sum _{i,j=1}^g \sum _{m= 1}^\infty   h^1(\mfr{x}_j,c(i,j,m ,\frac{\vare _i^m}{m} )/ t_j ,\mcal{O}_{\mfr{C}} )  \]
Through the isomorphism $\mrm{Im}(\partial _D)\cong \Lambda ^g $ given by the basis $\{ h^1(\mfr{x}_i,1/t_i ,\mcal{O}_{\mfr{C}})\} _i$ the composite
$(\log \circ \Phi )\colon \prod _{i=1}^g\mfr{m}^n \to \Lambda ^g$
is written as a power series map
\[ (\vare _i)_i\mapsto (\sum _{i=1}^g\sum _{m= 1}^\infty c(i,j,m)\frac{\vare _i^m}{m})_j . \]
Because of the tautological identity $c(i,j,1)=\delta _{ij}$, 
it induces the identity map on the tangent spaces at the origins;
thus it gives an analytic isomorphism around the basepoints of both spaces.
Since the taget $\Lambda ^g$ is homeomorphic to $\mrm{Pic}(\mfr{C})^{(n)}$ with its filtration topology (because the two are related by the log/exp isomorphism), the comparison of topologies is done.

\begin{ack}
Part of these notes was written as the author's 2013 master thesis.
During this work, the author was supported by Graduate Program for Leading Graduate Schools, MEXT, Japan.
When I started master's course in April 2011, there used to be frequent aftershocks of the earthquake in Tokyo too.
I thank Prof.\ Shuji Saito for guiding me firmly despite such circumstances. 
I am also grateful to Prof.\ Takao Yamazaki for encouragements to polish them up for publication.

\end{ack}

\begin{flushright}
\footnotesize {\it \today }
\normalsize
\end{flushright}


\begin{thebibliography}{EGA}


\bibitem[EGA]{EGA}
A. Grothendieck, J. Dieudonn\'{e}:
\'{E}l\'{e}ments de G\'{e}om\'{e}trie Alg\'{e}brique I, III premier/second parties.
Publ.\ Math.\ I.H.E.S., nos.\ 4, 11, 17. 1960, 1961, 1963.


\bibitem{Mattuck}
A. Mattuck:
Abelian varieties over $p$-adic ground fields.
Ann.\ of Math.\ vol.\ 62, no.\ 1, 
pp.92--119,
1955.


\bibitem{Milne}
J.S. Milne:
Jacobian Varieties.
In {\it Arithmetic Geometry}
(Storrs, Conn., 1984),
pp.167-212.
Springer-Verlag, New York, 1986.


\bibitem{Raynaud}
M. Raynaud:
Sp\'{e}cialisation du foncteur de Picard.
Publ.\ Math.\ I.H.E.S., no.\ 38,
pp.27--76, 1970.



\end{thebibliography}
\end{document}